\documentclass[a4paper,10pt]{article}
\def \Z {{\mathbf {Z}}}

\def\u{\bigsqcup}

\title{ \bf \Large  Chacon's Type Ergodic Transformations \\ with Unbounded Arithmetic Spacers. }
\author{V.V. Ryzhikov}
\date{E-mail: vryzh@mail.ru}

\begin{document}
\Large
\maketitle

{\large
 The following  generalizations of the Chacon map are proposed: instead of classical constant spacer sequence $(0,1,0)$ let  a sequence
$(0,s_j,0)$ be one with unbounded $s_j$. (We mention also an analogue of the  historical Chacon map with spacer sequences
in  the  form $(0,s_j)$.)  This narrow class of rank-one transformations 
is an  abundant  source of exercises,  open questions and maybe hard problems. All such constructions have partial rigidity, but some  other  properties could be different. For root sequence, $ s_j= [\sqrt{j}]$,
(or $ s_j= [\ln{j}]$)  the corresponding action is rigid, moreover it possesses all polynomials  in its weak closure.  In the linear case $s_j={j}$ we 
get (as  well as  for the classical Chacon transformation) the property of minimal  self-joinings (MSJ) which contrasts with the rigidity.  We present some observations   about
MSJ, mild mixing, partial mixing, $\ae$-mixing, absence  of factors, triviality of centralizer and spectral  primality.  We  state several problems, mention exponential  "self-similar" Chacon transformations  and flows on infinite measure spaces, and 
provide a place for the subsequent generalizations and the imagination. }

\section{Chacon-type maps} 
 { A rank-one construction} is determined by an integer $h_1$, a cut-sequence $r_j$ and  a spacer  sequence $\bar s_j$:
$$ \bar s_j=(s_j(1), s_j(2),\dots, s_j(r_j-1 ),s_j(r_j)), \ \ s_j(i)\geq 0.$$
We don't recall a  general definition,  below we  deal with  very
particular cases of rank-one constructions, i.e. 
 Chacon-type transformations, setting $h_1=1$,  $r_j=3$,  and  $$ \bar s_j=(0, s_j,0).$$  Let $s_j=1$ for all $j$,  the corresponding transformation is called Chacon's map. It has minimal self-joinings, so its  centralizer is trivial.  Chacon's map  is not rigid, but partially rigid and mildly  mixing \cite{JRS} (see also, \cite{PR}, \cite{JPRR}).    
\\
Let a sequence $s_j$ be fixed.
A Chacon-type transformation  $T$ on the step $j$ is defined as  a shift on a collection of disjoint sets (intervals)
$$E_j, TE_j, T^2E_j,\dots, T^{h_j}E_j.$$
\bf "Cutting" \rm $E_j$ into $3$ sets (subintervals)  of the same measure
$E_j=E_j^1\u E_j^2\u  E_j^3,$ 
 for  $i=1,2,3$ we  consider  columns
$E_j^i, TE_j^i ,T^2 E_j^i,\dots, T^{h_j}E_j^i.$
%
Adding $s_j$ so-called spacers over the second column we obtain 

$E_j^1, TE_j^1 ,T^2 E_j^1,\dots, T^{h_j}E_j^1,$

$E_j^2, TE_j^2 ,T^2 E_j^2,\dots, T^{h_j}E_j^2, T^{h_j+1}E_j^2,\dots T^{h_j+s_j}E_j^2,$

$E_j^3, TE_j^3 ,T^2 E_j^3,\dots, T^{h_j}E_j^3.$
\\
(all intervals   are disjoint and of the same measure).
We build $(j+1)$-tower  via  \bf "stacking" \rm :
$TT^{h_j}E_j^1 = E_j^{2}, \ \ \ \
TT^{h_j+s_j}E_j^2 = E_j^{3}.$ 
Now our transformation acts  on the $(j+1)$-tower
$$E_{j+1}, TE_{j+1} , T^2 E_{j+1},\dots, T^{h_{j+1}}E_{j+1},$$
where 
 $$E_{j+1}= E^1_j,\ \ \ \ h_{j+1}+1=3(h_j+1) +s_j.$$
  A limit cutting-and-stacking construction is an   invertible measure-preserving  map $T:X\to X$,
where $X$ is the union of all above intervals.
\\
The first thing that comes to mind is to consider:

1.  Root spacers, $ \bar s_j=(0, [\sqrt{j}],0);$

2.  Linear spacers, $ \bar s_j=(0, {j},0);$

3. Polynomial spacers, $ \bar s_j=(0, {j^m},0);$

4.  Exponential spacers, $ \bar s_j=(0, {2^j},0);$

5. Prime  spacers,  $ \bar s_j=(0, {p(j)},0),$
where $p(j)$  denotes $j$-th prime number.
\\
The second  thing  is to produce  \bf 
systems with   "spacer memory": \rm
$$ \bar s_{j_k}=(0, h_{m_k},0).$$
Such transformations could be rigid even as  $ s_j\to +\infty$ and $(s_{j+1} -s_j)\to +\infty$.
In fact root Chacon map has a nice spacer memory. 
We will not discuss this intriguing phenomenon, let us 
 note only  that the transformation properties "vary" considerably,
if the intersection of sets $\{s_j\}$  and $\{h_j\}$  is infinite.

\section{ Properties in question} 
Below, we recall some  asymptotic, spectral and joining properties of transformations, which will be of interest. 
\\
\bf MSJ. \rm A {\it self-joining} (of order 2) is defined to be a $T\times T$-invariant
measure $\nu$ on $X\times X$ with the marginals  equal to  $\mu$:
$$\nu(A\times X)=\nu(X\times A)=\mu(A).$$
A joining $\nu$ is called ergodic if the dynamical system
$(T\times T, X\times X, \nu)$ is ergodic.
The measures $\Delta^i= (Id\times T^i)\Delta$ 
defined by the formula
$
	\Delta^i(A\times B) = \mu(A\cap T^iB). 
$ 
If $T$ is ergodic, then $\Delta^i$
is an   ergodic self-joining. 

We say that  $T$  has {\bf minimal self-joinings}
  ($T$ has  MSJ)  if the set of all  its ergodic
self-joinings contains only  $\mu\times\mu$ and $\Delta^i$, $i\in \Z$.

The property of minimal self-joinings implies {  mild mixing}, and  { partial rigidity} for non-mixing transformations. We recall that an automorphism $T$ is {\bf mildly mixing} if  for any set $A$, $0<\mu(A)<1$,  \
$\limsup_j 
\mu(A\cap T^iA)  < \mu(A) .$
An automorphism  $T$ is mildly mixing iff it has no rigid factors ($S$ is rigid, if there is a sequence $m_j\to\infty$
such that $S^{k_j}\to Id$).  

 An automorphism $T$ is  {\bf partially rigid} with a coefficient $\rho$, if
 there is a sequence $n(i)\to\infty$ such that  for any  measurable sets $A,B$ \
$\lim_{i\to\infty} 
\mu(A\cap T^{n(i)}A)  \geq \rho \,\mu(A).$
The maximum of such   $\rho$-s is denoted by $\rho(T)$.
If $\rho(T)=1$, then $T$ is called \bf rigid.\rm

A map $T$ is  \bf  $\ae$-mixing, \rm  $\ae\in (0,1)$,  if there is a sequence  $k_j$ such that
in $L_2^0(X,\mu)$   the operators $T^{k_j}$ weakly  converge to $(1-\ae)I$ (\cite{O},\cite{S}).

If all symmetric powers of $T$ have simple spectrum, we say that
$T$ has \bf prime  spectrum.  \rm This property, as well as  $\ae$-mixing property, implies \bf DC, \rm the pairwise disjointness of the convolutions $\sigma^{\ast n}$, $n>1$, for  spectral measure $\sigma$ of $T$.  

Classical Chacon's  transformation is mildly mixing (it has MSJ \cite{JRS}), is not $\ae$-mixing \cite{JPRR},  DC of Chacon's  transformation has been proved by A. Prikhodko and the author; F. Parreau,  O. Ageev established the primality  of its spectrum.

\vspace{3mm}
\section{Observations, conjectures, questions} 
\it
1.  Root Chacon map is rigid and mildly mixing.  It possesses all
polynomial weak limits of powers. So it has prime spectrum. 
\\
2.  Linear Chacon map  has MSJ (weakly mixing,
non-rigid), $\ae$-mixing,  $\rho(T)\geq \frac 3 4$. 
\\
3. Polynomial Chacon map is  mildly mixing. (Trivial centralizer, no factors.)
Conjecture: $T$ has MSJ.    Prime spectrum?
\\
4.  Exponential Chacon map.   ???
\\
5. Prime  Chacon map.   It is  $\ae$-mixing (prime spectrum is very plausible).   Is it partially mixing?
Conjecture:  MSJ.  \rm

\bf 2-adic Chacon's transformations.  \rm Similar problems appear in case of  $r_j=2$,  $\bar s_j=(0, s_j))$.

\it
1.  Root  2-Chacon map $( \bar s_j=(0, [\sqrt{j}]))$ is $\ae$-mixing.
   Is it non-rigid? Conjecture:  prime spectrum.

2.  Linear 2-Chacon map $( \bar s_j=(0, {j}))$  is $\ae$-mixing and mildly mixing.
(Trivial centralizer, no factors.)
  Conjecture: MSJ, prime spectrum.

3. Polynomial 2-Chacon map $ (\bar s_j=(0, {j^m}), m\geq 2)$ mildly mixing. (Trivial centralizer, no factors.)
Conjecture:  MSJ.

4.  Exponential 2-Chacon map $( \bar s_j=(0, {2^j}))$.  Now, the corresponding  transformation acts on an \bf infinite measure \rm space. (To have a finite measure the reader could modify the construction.) ???  

4. Prime 2-Chacon map,  $( \bar s_j=(0, {p(j)}))$. It is weakly mixing.
Questions: partial mixing (MSJ), $\ae$-mixing, prime spectrum?
\rm

{\bf Around Chacon-type transformations.}
  Rank-one  constructions  are naturally divided by 4  classes:   Unbounded cutting + Unbounded spacer, Bounded cutting + Unbounded spacer, ...
In connections with the work \cite{B}, where   BB-class  ( Bounded cutting + Bounded spacer) 
was considered, we proved \cite{R12}  that  mildly mixing BB-constructions
have MSJ property.  What about the BU-constructions? Is it true that  (positive) powers of weakly mixing BU-construction are pairwise disjoint?

\bf Conjecture. \it
Weakly mixing BB-constructions have prime spectrum, but  never $\ae$-mixing. \rm \\
King and Thouvenot \cite{KT} proved MSJ for partially mixing
rank-one maps. 

\bf Conjecture:  \it there is a partially mixing unbounded  Chacon transformation.\rm
\\
(Recall that all rank-one constructions with bounded spacer  are  not partially mixing.)

\section{ Root Chacon transformation} 
\it Root Chacon map is rigid and weakly mixing.  It possesses all
polynomial weak limits of powers.\rm

Proof. First, we show the rigidity.
Let $\rho(T)=\rho$, we have $$T^{n(k)} \to \rho I +\dots, \ \ \ n(k)\to +\infty.$$
On large intervals our spacer sequence $s_j=[\sqrt {j}]$  is constant.
Let 
$$ [\sqrt {j(k)}]=[\sqrt {j(k)+1}]=[\sqrt {j(k)+2}]=\dots =[\sqrt {j(k)+m(k)}], \ m(k)\to\infty,$$
then   $$T^{-h_j} \approx_w \frac 1 2  (I + T^{[\sqrt j]}).$$
Setting $[\sqrt {j(k)}]= n(k)$, we get
$$T^{-h_{j(k)}} \approx_w \frac 1 2  (I + T^{n(k)})
\approx_w \frac 1 2  I + \frac 1 2 \rho I  +\dots .$$
Thus,
$$ \rho(T)\geq  \frac 1 2  +\frac 1 2  \rho(T), \ \ \ \rho(T)=1,$$
 our $T$ is rigid.  Let
 $$T^{r(k)} \to  I.$$ Given $p$
we take  $[\sqrt {j(k)}]= r(k)+p$.  Moreover,  let  
$$ r(k)+p =[\sqrt {j(k)}]=[\sqrt {j(k)+1}]=[\sqrt {j(k)+2}]\dots =[\sqrt {j(k)+m(k)}], \ m(k)\to\infty.$$
Then
$$ T^{-h_{j(k)}} \approx_w  \frac 1 2 I +\frac 1 2  T^p.  $$
Thus, for any $p$ we find $q(k)$ such that
$$ T^{q(k)} \to  \frac 1 2 I +\frac 1 2  T^p.  $$
The  latter implies  that the weak closure of $\{T^n\}$
contains all polynomials \cite{R}. 

 We get the same  by $s_j=[{j}^{1-\varepsilon}]$, $s_j=[\log{j}]$, $s_j=[\log\log{j}]$ and so on.
\\ \\ \\ 
{\bf Assertion. }{ \it Let    intervals(!)    $I_s =\{j: s_j=s\}$  satisfy  $|I_s|\to\infty$.   Then the corresponding Chacon transformation possesses all
polynomial weak limits of powers and, hence, it has prime spectrum.\rm
}
\section{ Linear Chacon map} 
\bf Linear Chacon map  $T$ is not partially mixing. \rm It is quasi-rigid: for any $\varepsilon >0$
there is a weak limit in the form    $\sum_{i\geq 0} \alpha_iT^i +\dots$, where $\sum_{i\geq 0} \alpha_i > 1-\varepsilon $.

To see this let us consider  $$\alpha(T) =\sup\{\sum_{i\geq 0} \alpha_i \ :\ \left(\sum_{i\geq 0} \alpha_iT^i +\dots\right)\in Lim(T)  \}.$$  Reasoning as   above,  we get
$$\alpha(T)\geq  \frac 1 2  +\frac 1 2  \alpha(T),    \ \ \  \alpha =1.$$
\\
\bf Linear Chacon map   is $\ae$-mixing. \rm 
It's not hard to see that
$T$ is weakly mixing. There is a sequence $j(k)\to\infty$  such that
$$ T^{h_{j(k)}}\to \frac 1 2 I +\frac 1 2\Theta$$ 
weakly  in  $L_2(X,\mu),$
where $\Theta f\equiv \int f d\mu$.
\\ \\
\bf Linear Chacon map   is  mildly mixing. \rm 
Let $T^{n(i)}\to P$. We have the following partial factorization: 
 $$P=\frac Q {3(3I-T)} + \frac 5 6 R=\frac 1 9  Q \left(I +\frac T 3+\frac {T^2} 9+\dots\right)+  \frac 5 6 R \eqno (I,T),$$
where  $Q,R$are some Markov operators. 
 
From King's result we know that any proper factor of a rank-one map have to be rigid. So there  is a sequence $n(i)\to \infty$ for which $T^{n(i)}|F\to I$ in the Hilbert space $F$ corresponding to  the  factor. Let $T^{n(i')}\to P$. We have 
$ I|F= \frac 1 9  Q \left(I +\frac T 3\right)|F+\dots$.  But this is possible only in the case of trivial factor, where $T|F=I|F$. 
From $(I,T)$ and King's weak closure theorem  it follows the triviality of the centralizer of $T$.  Thus, $T$ is mildly mixing.\\
\\
\bf Linear Chacon map  has MSJ. \rm Use a  modification of the proof of MSJ for mildly mixing BB-constructions \cite{R12}.

 \section{Polynomial Chacon map} It is mildly mixing (thus, its centralizer and factors are all trivial) and $\ae$-mixing.   If any  limit of its powers has the form $aT^n+(1-a)\Theta$, then polynomial Chacon map is partially mixing and has MSJ.

\section {Prime  Chacon map} It is totally ergodic. Due to  Y. Zhang \cite{Z} we easily get 
the weak mixing via limits in the form  
$$ \frac 1 3 I+\frac 1 9 T^{d}+ \frac 5 9  P, \ \ d>1.$$ 

From the recent work by
 Maynard \cite{M} it follows for any $m> 1$
the existence of  weak limits 
$$ \frac 1 3 I+\frac 1 9 T^{d_1}+\dots +\left(\frac 1 3 \right)^{m+1}T^{d_m}+ b P_m$$ for some Markov operators $P_m$ and $d_m > d_{m-1}>\dots>d_1>1$.

\vspace{5mm}
Prime  Chacon map is $\ae$-mixing (so it has DC):   for  some  sequence $\{j(k)\}$
$$ T^{h_{j(k)}}\to \frac 1 2 I +\frac 1 2\Theta$$ 
(because of the fact that for   a weakly mixing $T$ the  powers $T^{p}$ are close to $\Theta$  for most primes $p$).
\\ Questions.  Is it true that 
$ T^{h_{j}}\to \frac 1 2 I +\frac 1 2\Theta \ \ as \ \  j\to\infty \ \ ?$ Is Prime  Chacon map partially mixing?
 \\ \\
 Conjecture. \it Prime  Chacon map possesses a weak limit in the form
$$ \frac 1 3  (I+\frac 1 3 T^{d_1}+\dots +\frac  1 {3^m} T^{d_m}+\dots) + \frac 1 2 \Theta,$$
and it  has prime spectrum.\rm
\newpage
\section{Exponential  "self-similar" Chacon maps} Given  $n>3$ let $h_1=1$, $\bar s_j=(0,s_j,0)$, where  $s_1=1$ and $s_j=n^j-3n^{j-1}$ as $j>1$.  The corresponding
 transformation $T$ preserves infinite measure, it has simple continuous self-similar spectrum.
\\
\bf Assertion 1.  \it The power $T^{n^d}$ has homogeneous spectrum of multiplicity ${n^d}$. 

 \rm Proof.  Indeed, 
$T^n\cong T\oplus \dots\oplus T\ (n\ times)$. \rm
\\
To get the same for $n=3$ we apply 2-adic Chacon map. \\
\bf Assertion 2.  \it If $n>2$,  $h_1=1$, $s_1=1$ and  $\bar s_j=(0,n^j-2n^{j-1})$ as $j>1$, then 
$$T^n\cong T\oplus \dots\oplus T\ (n\ times).$$ \rm

Quite similarly, we obtain self-similar rank-one  flows ($T_t\cong T_{nt}$), acting on spaces with infinite measure. The associated Poisson suspensions give some new  examples of "finite" self-similar flows.

\end{document}